\documentclass{birkjour}

\usepackage{amsmath,amsthm}
\usepackage[english]{babel}
\usepackage{amsfonts}
\usepackage{latexsym}
\usepackage{times}
\usepackage{tikz-cd}

\newcommand{\C}{{\mathbb C}}
\newcommand{\lb}{\lambda}
\newcommand{\bb}{{\mathcal B}}

\newcommand{\pf}{{\bf Proof. }}

\newcommand{{\rr}}{\mathcal{R}}
\newcommand{\la}{\langle}
\newcommand{\ra}{\rangle}

\newtheorem{thm}{Theorem}[section]
\newtheorem{corr}[thm]{Corollary}
\newtheorem{lem}[thm]{Lemma}
\newtheorem{prop}[thm]{Proposition}
\newtheorem{exam}[thm]{Example}

\newcommand{\zb}{$\hfill\Box$}

\renewcommand{\theequation}{\thesection.\arabic{equation}}
\makeatletter\@addtoreset{equation}{section} \makeatother

\begin{document}

\title {Hermitian geometry on resolvent set (I)}

\author[R. G. Douglas]{Ronald G. Douglas}
\address{Department of Mathematics, Texas A\&M University\\ College Station, TX 77843, U.S.A.} \email{rdouglas@math.tamu.edu}

\author[R. Yang]{Rongwei Yang}
\address{School of Mathematical Sciences,
Tianjin  Normal University\\ Tianjin 300387, P.R. China.\\ AND \\
Department of Mathematics and Statistics, SUNY at Albany\\
Albany, NY 12222, U.S.A.} \email{ryang@albany.edu}



\begin{abstract}
For a tuple $A=(A_1,\ A_2,\ ...,\ A_n)$ of elements in a unital  Banach algebra ${\mathcal B}$, its {\em projective joint spectrum}
$P(A)$ is the collection of $z\in \C^n$ such that $A(z)=z_1A_1+z_2A_2+\cdots +z_nA_n$ is not
invertible. It is known that the ${\mathcal B}$-valued $1$-form $\omega_A(z)=A^{-1}(z)dA(z)$ contains much topological information about the joint resolvent set $P^c(A)$. This paper studies geometric properties of $P^c(A)$ with respect to Hermitian metrics defined through the ${\mathcal B}$-valued {\em fundamental form} $\Omega_A=-\omega^*_A\wedge \omega_A$ and its coupling with {\em faithful} states $\phi$ on $\bb$, i.e. $\phi(\Omega_A)$. The connection between the tuple $A$ and the metric is the main subject of this paper. In particular, it shows that the K\"{a}hlerness of the metric is tied with the commutativity of the tuple, and its completeness is related to the Fuglede-Kadison determinant.
\end{abstract}

\maketitle

\section{Introduction}

In \cite {CD}, Cowen and the first-named author introduced geometric concepts such as holomorphic bundle and curvature into Operator Theory. This gave rise to complete and computable invariants for certain operators with rich spectral structure, i.e, the Cowen-Douglas operators. The idea was extended to multivariable cases in the study of Hilbert modules in analytic function spaces, where curvature invariant is defined for some canonical tuples of commuting operators. We refer readers to \cite{DM,Do,DP, Sa} and the references therein for more information. This paper aims to study general (possibly non-commuting) tuples using geometric ideas. The approach is based on the newly emerged concept of projective joint spectrum $P(A)$ introduced by the second-named author (cf. \cite{Ya}). Some preliminary findings are reported here. They shall help to build a foundation for subsequent and more in-depth studies.

 Let ${\mathcal B}$ be a complex Banach algebra with unit $I$ and $A=(A_1, A_2, \cdots, A_n)$ be a tuple of linearly independent elements in ${\mathcal B}$. The multiparameter pencil \[A(z):=z_1A_1+z_2A_2+\cdots +z_nA_n\] is an important subject of study in various fields, for example in algebraic geometry \cite{Vi}, math physics \cite{AJ,Je}, PDE \cite{At}, group theory \cite{GS}, etc., and more recently in the settlement of the Kadison-Singer Conjecture (\cite{MSS}). In these studies, the primary interest was in the case when $A$ is a tuple of self-adjoint operators. For general tuples, the notion of projective joint spectrum is defined in \cite{Ya} and its property is investigated in a series of papers (cf. \cite{BCY}, \cite{CSZ}, \cite{CY}, \cite{SYZ}). \\

{\bf Definition}. For a tuple $A=(A_1, A_2, \cdots, A_n)$ of elements in a unital Banach algebra ${\mathcal B}$, its projective joint spectrum is defined as \[P(A)=\{z\in {\C}^n| A(z)\ \text{is not invertible}\}.\] The {\em projective resolvent set} refers to the complement $P^c(A)={\C}^n\setminus P(A)$. \\

Various notions of joint spectra for (commuting) tuples have been defined in the past, for instance the Harte spectrum (\cite{Ha72,Ha73}) and the Taylor spectrum (\cite{Fa,Ta70,Ta72}), and they are key ingredients in multivariable operator theory. The projective joint spectrum, however, has one notable distinction: it is ``base free" in the sense that, instead of using $I$ as a base point and looking at the invertibility of $(A_1-z_1I,\ A_2-z_2I,\ ...,\ A_n-z_nI)$ in various constructions, it uses the much simpler pencil $A(z)$. This feature makes the projective joint spectrum computable in many interesting noncommuting examples, for example the tuple of general compact operators (\cite{SYZ}), the generating tuple for the free group von Neumann algebra (\cite{BCY}), the Cuntz tuple of isometries $(S_1,\ S_2,\ \dots,\ S_n)$ (\cite{Cu,HWY}) satisfying
\begin{eqnarray*}
&&S_i^*S_j=\delta_{ij}I \quad \textrm{for}\ 1\leq i, j\leq n,\\
&&\sum_{i=1}^{n}S_iS_i^*=I,
\end{eqnarray*}
where $I$ is the identity. For the generating tuple $(1,\ a,\ t)$ of the infinite dihedral group \[ D_{\infty}=<a,\ t\ |\ a^2=t^2=1>,\] the projective spectrum has found convincing applications to the study of group of intermediate growth (finitely generated group whose growth is faster than polynomial growth but slower than exponential growth (cf. \cite{GY})). 

For a general tuple $A$, it can happen that $P(A)=\C^n$ and hence $P^c(A)$ is empty. But in this paper we shall always assume that $P^c(A)$ is non-empty. It is known that every path-connected component of $P^c(A)$ is a domain of holomorphy. This fact is a consequence of a result in \cite{ZKKP}, and it is also proved independently in \cite{HY}. From this point of view, $P^c(A)$ has good analytic properties. On $P^c(A)$, the holomorphic Maurer-Cartan type $\bb$-valued $(1,0)$-form\[\omega_A(z)=A^{-1}(z)dA(z)=\sum_{j=1}^nA^{-1}(z)A_jdz_j\] 
played an important role in \cite{Ya}, where it is shown to contain much information about the topology of $P^c(A)$. For example, the de Rham cohomology $H^*(P^c(A),\C)$ can be computed by coupling $\omega_A$ with invariant multilinear functionals or cyclic cocycles (\cite{CY,Ya}). Two simple facts about $\omega_A$ are useful here. 

1) By differentiating the equation $A(z)A(z)^{-1}=I$, one easily verifies that 
\[dA^{-1}(z)=-A^{-1}(z)(dA(z))A^{-1}(z),\]
and consequently,
\begin{equation}
d\omega_A=dA^{-1}(z)\wedge dA(z)=-\omega_A\wedge \omega_A.
\end{equation}
Here $d$ is complex differentiation (see Section 1) and $\wedge$ is the wedge product. 

2) $\omega_A$ is invariant under left multiplication on the tuple $A$ by invertible elements. To be precise, we let $GL(\bb)$ denote the set of invertible elements in $\bb$. Then for any $L\in GL(\bb)$ and $B=LA=(LA_1,\ LA_2,\ \cdots,\ LA_n)$, we have for the multiparameter pencils $B(z)=LA(z)$, and hence $P(B)=P(A)$. Moreover, 
\begin{equation}
\omega_B=B^{-1}(z)dB(z)=A^{-1}(z)L^{-1}LdA(z)=\omega_A.
\end{equation}
This shows that $\omega_A$ is invariant under the left action of $GL(\bb)$ on the tuple $A$. 

For simplicity, we shall use Einstein convention for summation in many places in this paper. For example, we shall write
$\omega_A(z)=A^{-1}(z)A_jdz_j.$ When ${\mathcal B}$ is a $C^*$-algebra, the adjoint of $\omega_A(z)$ is the $\bb$-valued $(0,1)$-form
\[\omega^*_A(z)=(A^{-1}(z)A_j)^*d\overline{z}_j.\] 
Now we define the {\em fundamental form} for the tuple $A$ as 
\begin{equation}
\Omega_A=-\omega^*_A\wedge \omega_A=(A^{-1}(z)A_j)^*A^{-1}(z)A_kdz_k\wedge d\overline{z}_j.
\end{equation}
Here the negative sign exists in the first equality because it is conventional to write a $(1,\ 1)$-form as linear combinations of $dz_j\wedge d\bar{z}_k$. Often, the factor $\frac{i}{2}$ is used to make $\frac{i}{2}\Omega_A$ a {\em self-adjoint} form (cf. (1.3)), but it is not important in this paper.

For a suitable choice of linear functional $\phi$, such as a faithful state on ${\mathcal B}$, the evaluation \[\phi(\Omega_A)=\phi\left( (A^{-1}(z)A_j)^*A^{-1}(z)A_k\right)dz_k\wedge d\overline{z}_j\] induces a positive definite bilinear form on the holomorphic tangent bundle of $P^c(A)$, thus giving a Hermitian metric on $P^c(A)$. The connection between the metric and the tuple $A$ shall be the primary concern of this paper. In particular, it shows that the K\"{a}hlerness of the metric is tied with the commutativity of the tuple (cf. Theorem 2.4). A notable feature of this metric is that it has singularities at the joint spectrum $P(A)$. So completeness of the metric is an important issue. 

\section{$\bb$-valued differential forms}

Let $M$ be a complex manifold of dimension $n$. If $z=(z_1, z_2, \cdots, z_n)$ is the coordinate in a local chart, then $\partial_i$ stands for $\frac{\partial}{\partial z_i}$, and $\bar{\partial}_i$ stands for $\frac{\partial}{\partial \bar{z}_i}$. As a convention, we let $\partial=\sum_i \partial_i$, $\bar{\partial}=\sum_i \bar{\partial}_i$, and $d=\partial+\bar{\partial}$. Consider a globally defined smooth $(1, 1)$-form $\Phi(z)=g_{jk}(z)dz_j\wedge d\bar{z_k}$ expressed in the local chart with the Einstein summation convention. $\Phi$ induces a bilinear form $\hat{\Phi}(z)=g_{jk}(z)dz_j\otimes d\bar{z_k}$ on the holomorphic tangent bundle $T(M)$ over the local chart such that 
\[\hat{\Phi}(z)(\partial_j,\ \partial_k)=g_{jk}(z).\]
We say that $\Phi$ defines a Hermitian metric on $M$ if the $n\times n$ matrix function $g(z)=(g_{jk}(z))$ is positive definite for each $z\in M$. In this case $g(z)$ is called the associated metric matrix, and $\frac{i}{2}\Phi(z)$ is called the fundamental form of the metric.  Here, the constant $\frac{i}{2}$ is to normalize the fundamental form so that \\

\hspace{2mm} (a) $\frac{i}{2}\Phi(z)$ is real i.e. $\frac{i}{2}\Phi(z)=\overline{\frac{i}{2}\Phi(z)}$;

\hspace{2mm} (b) in the one variable case $z=x+iy$, we have $\frac{i}{2}dz\wedge d\bar{z}=dx\wedge dy$. \\

Some geometric concepts are relevant to the study here.\\

\hspace{2mm} 1) A Hermitian metric induced by $\Phi(z)$ is said to be K\"{a}hler if $d\Phi=0$, i.e. $\Phi(z)$ is closed.

\hspace{2mm} 2) The Ricci curvature tensor is the $n\times n$ matrix $R(z)=(R_{jk}(z))$, where $R_{jk}=-\partial_j\bar{\partial_k} \log \det g(z)$, and the Ricci form is $Ric_{\Phi}(z)=R_{jk}(z)dz_j\wedge d\bar{z_k}$. 

\hspace{2mm} 3) The metric is said to be Ricci flat if $R(z)=0$ for each $z\in M$. 

\hspace{2mm} 4) The metric is said to be Einstein if $Ric_{\Phi}=\lambda \Phi$ for some constant $\lambda$, and it is Calabi-Yau if it is K\"{a}hler and Ricci flat (i.e. $\lb=0$).\\

One checks that a Hermitian $(1, 1)$-form $\Phi(z)=g_{ij}(z)dz_i\wedge d\bar{z_j}$ induces a K\"{a}hler metric if and only if \begin{equation}
\partial_kg_{ij}(z)=\partial_ig_{kj}(z),\ \ z\in M,
\end{equation}
 for all $1\leq i, j, k\leq n$. Further, $\Phi$ is Ricci flat if and only if $\log \det g(z)$ is pluri-harmonic (cf. \cite{Ra}). Hence in this case $\log \det g(z)=f(z)+\overline{f(z)}$ holds locally for some holomorphic function $f$, i.e., 
\begin{equation}
\det g(z)=e^{f(z)+\overline{f(z)}}=|e^{f(z)}|^2.
\end{equation}

Now we consider a smooth $\bb$-valued $(p,q)$-form $\omega^{p,q}(z)$ defined on a complex domain $R\subset \C^n$, where $0\leq p,q\leq n$. The set of such forms is denoted by $\Lambda^{p,q}(R, {\mathcal B})$. Then $\omega^{p,q}(z)$ is said to be $\partial$-closed or $\overline{\partial}$-closed if $\partial \omega^{p,q}(z)=0$, or respectively  $\overline{\partial} \omega^{p,q}(z)=0$ on $R$. If $\omega^{p,q}$ is both $\partial$-closed and $\overline{\partial}$ closed, then $d\omega^{p,q}=0$, and in this case we simply say $\omega^{p,q}$ is closed.

For the rest of the paper, unless stated otherwise, we shall assume ${\mathcal B}$ is a $C^*$-algebra. Recall that for a tuple $A$, the fundamental form $\Omega_A$ is given by (0.3). Now we check that $\frac{i}{2}\Omega_A$ is {\em self-adjoint} in the sense that 
\begin{equation}
(\frac{i}{2}\Omega_A)^*=\frac{i}{2}\Omega_A.
\end{equation}
 Indeed, one checks by (0.3) that
\begin{align*}
(\frac{i}{2}\Omega_A)^*&=\frac{-i}{2}(-\omega^*_A\wedge \omega_A)^{*}\\
&=\frac{-i}{2}\left((A^{-1}(z)A_j)^*A^{-1}(z)A_kdz_k\wedge d\overline{z}_j\right)^*\\
&=\frac{-i}{2}\left(A^{-1}(z)A_k)^*A^{-1}(z)A_jd\overline{z}_k\wedge dz_j\right)\\
&=\frac{-i}{2}(\omega^*_A\wedge \omega_A)=\frac{i}{2}\Omega_A. 
\end{align*}
Further, in the case $n=1$ and $A_1$ is invertible, one has $A(z)=zA_1$, $P^c(A)=\C\setminus \{0\}$, and $\omega_A(z)=\frac{dz}{z}$. Writing $z=x+yi$ we have
\[\frac{i}{2}\Omega_A=\frac{i(dz\wedge d\overline{z})}{2|z|^2}=\frac{dx\wedge dy}{x^2+y^2},\]
which is real. For a linear functional $\phi\in {\mathcal B}^*$,
\[\phi(\frac{i}{2}\Omega_A)=\frac{i}{2}\phi\left(A^{-1}(z)A_j)^*A^{-1}(z)A_k\right)dz_k\wedge d\overline{z}_j,\]
which is a regular $(1,1)$-form. Moreover, if $\phi$ is positive then by (1.3)
\[\overline{\phi(\frac{i}{2}\Omega_A)}=\phi(\frac{-i}{2}\Omega^*_A)=\phi(\frac{i}{2}\Omega_A),\]
and hence $\phi(\frac{i}{2}\Omega_A)$ is a regular real $(1,1)$-form. The following definitions are ${\bb}$-valued versions of Hermitian metric and K\"{a}hler metric. \\ 

{\bf Definition}. A smooth $\bb$-valued $(1,1)$-form $\Phi(z)=\Phi_{jk}(z)dz_j\wedge d\overline{z}_k$ on a complex domain $R$ is said to be Hermitian if for every $z\in R$ and any set of $n$ complex numbers $v_1,\ v_2,\ \dots,\ v_n$ the double sum $v_i\Phi_{jk}(z)\overline{v_k}$ is a positive element in $\bb$.\\

{\bf Definition}. A smooth $\bb$-valued Hermitian $(1,1)$-form $\Phi(z)$ is said to be K\"{a}hler if it is closed, i.e. $d\Phi(z)=0$.\\

Two nuances are worth mentioning. First, there is a difference between a self-adjoint $\bb$-valued $(1,1)$-form as defined in (1.3) and a Hermitian $\bb$-valued $(1,1)$-form defined above--the latter requires positivity of $\bb$-valued matrix $(\Phi_{ij}(z))$. Second, in the scalar-valued case, for a Hermitian metric defined by $\Omega(z)=g_{jk}(z)dz_j\wedge d\bar{z}_k$, the complex metric matrix $(g_{jk}(z))$ needs to be positive definite, but in the ${\bb}$-valued case, it is only required that $v_j\Phi_{jk}(z)\overline{v_k}\geq 0$ in $\bb$. 
For the fundamental form $\Omega_A$ and any set of $n$ complex numbers $v_1,\ v_2,\ \dots,\ v_n$, one checks by (0.3) that the double sum
\begin{align*}
v_jA_k^*(A^{-1}(z))^*A^{-1}(z)A_j\bar{v_k}&=\bar{v_k}A_k^*(A^{-1}(z))^*A^{-1}(z)v_jA_j)\\
&=(A^{-1}(z)A(v))^*(A^{-1}(z)A(v)\geq 0.\tag{1.4}
\end{align*}
So $\Omega_A $ is Hermitian. 

In this paper, the extended tuple $\hat{A}=(I,\ A_1,\ A_2,\ \cdots,\ A_n)$ is used in some places for normalization purposes. In this case we can regard $I$ as $A_0$ and $P(\hat{A})$ is then clearly a subset in $\C^{n+1}$. This extension guarantees two convenient facts: 

1) obviously, $(1,0,0,\dots,0)\in P^c(\hat{A})$;

2) the identity $I$ is in the range of $\hat{A}(z)$ as $z$ varies in $\C^{n+1}$. 

Now we are in position to state the first result.
\begin{prop}
For the tuple $A=(A_1,\ A_2,\ \cdots,\ A_n)$, the following statements are equivalent.

(i) $A$ is commuting.

(ii) $\omega_{A}(z)$ is closed.

(iii) $\Omega_{\hat{A}}(z)$ is K\"{a}hler.
\end{prop}
\pf Using (0.1) and the fact that $dz_j\wedge dz_k=-dz_k\wedge dz_j$, one checks that
\begin{align*}
d\omega_A&=-\omega_A \wedge \omega_A\\
&=A^{-1}(z)A_j\ A^{-1}(z)A_kdz_j\wedge dz_k\\
&=\sum_{j<k}[A^{-1}(z)A_j,\ A^{-1}(z)A_k]dz_j\wedge dz_k,\tag{1.5}
\end{align*} 
where $[a, b]=ab-ba$. So when $A$ is a commuting tuple we have 
\[[A^{-1}(z)A_j,\ A^{-1}(z)A_k]=0,\ \ z\in P^c(A),\]
and consequently $d \omega_A=\partial \omega_A=0$. This shows (i) implies (ii). 

Recall that $\Omega_A(z)=-\omega^*_{A}\wedge\omega_{A}$. Since $\omega_A(z)$ is holomorphic
\begin{align*}
d\Omega_A&=-(\partial +\bar{\partial})(\omega^*_A\wedge \omega_A)\\
&=-(\partial \omega_A^*+\bar{\partial}\omega^*_A)\wedge \omega_A+\omega_A^*\wedge (\partial \omega_A +\bar{\partial}\omega_A)\\
&=-\big((\partial \omega_A)^*\wedge \omega_A \big)+(\omega^*_A\wedge \partial \omega_A). \tag{1.6}
\end{align*}
So if $\omega_{A}(z)$ is closed then $d\Omega_A=-(0 \wedge \omega_A)+(\omega^*_A\wedge 0)=0$, i.e. $\Omega_A$ is K\"{a}hler. Since $\hat{A}$ is commuting, the forgoing arguments in particular implies that $\Omega_{\hat{A}}$ is K\"{a}hler as well.
This shows (ii) implies (iii). Now if $\Omega_{\hat{A}}$ is K\"{a}hler then by (1.6)
\begin{align*}
0=d\Omega_{\hat{A}}=-(\bar{\partial}\omega_{\hat{A}}^*)\wedge \omega_{\hat{A}}+\omega^*_{\hat{A}}\wedge\partial  \omega_{\hat{A}}.\tag{1.7}
\end{align*}
Since the first term in (1.7) is a $(1,2)$-form and the second is a $(2,1)$-form, the sum is $0$ only if both terms are $0$. Furthermore, by (1.5) and (1.7) we have
\begin{align*}
0&=\omega^*_{\hat{A}}\wedge d \omega_{\hat{A}}\\
&=(\hat{A}^{-1}(z)A_i)^*\sum_{0\leq j<k}[\hat{A}^{-1}(z)A_j,\ \hat{A}^{-1}(z)A_k]d\bar{z}_i\wedge dz_j\wedge dz_k, 
\end{align*}
which implies that on $P^c(\hat{A})$,
\begin{equation*}
(\hat{A}^{-1}(z)A_i)^*[\hat{A}^{-1}(z)A_j,\ \hat{A}^{-1}(z)A_k]=0,\ \ \ \forall \ 0\leq i,\ j,\ k\leq n.
\end{equation*}
Setting $i=0$ and using the fact $\hat{A}(z)=A_0=I$ at $z=(1,0,0,\dots,0)\in \C^{n+1}$, we have \[A_kA_j=A_jA_k,\ \forall 1\leq j,\ k\leq n,\] i.e. the tuple $A$ is commuting. So (iii) implies (i).
\zb

To prove that (iii) implies (i), we used the fact that the tuple $\hat{A}$ includes $I$. In fact, as the proof indicates that if $A(z)=I$ for some $z\in \C^n$ then there is no need to use $\hat{A}$.  But this kind of requirement is indispensible.

\begin{exam}
Consider a pair $A=(A_1,\ A_2)$, where $A_1$ is invertible. Then writting $A(z)=A_1(z_1+z_2A_1^{-1}A_2)$ and using (0.2), we have \[\omega_A=(z_1I+z_2A_1^{-1}A_2)^{-1}d(z_1I+z_2A_1^{-1}A_2).\]
Since $I$ commutes with $A_1^{-1}A_2$, the fundamental form $\Omega_A$ is K\"{a}hler by Proposition 1.1. But clearly $A_1$ may not commute with $A_2$.
\end{exam}
For every $\phi\in {\mathcal B}^*$, we have $d\phi(\omega_A)=\phi (d\omega_A)$. We let $H^q(P^c(A), \C)$ denote the $q$-th de Rham cohomology of $P^c(A)$. The following is then an easy consequence of Propositions 1.1.

\begin{corr}
Let $A$ be a commuting tuple. Then $\phi(\omega_A)\in H^1(P^c(A),\C)$ for every $\phi\in \bb^*$.
\end{corr}

Proposition 1.1 raises a natural question: is $\omega_{A}$ is exact i.e. is there is a $\bb$-valued smooth function $f(z)$ on $P^c(A)$ such that $\omega_{A}=df$? We take time to address this issue now. For every $\phi\in \bb^*$ and a $\bb$-valued smooth $(p, q)$-form $\omega^{p,q}(z)$ on a complex domain $R$, the coupling $\phi(\omega^{p,q}(z))$ is a scalar-valued smooth $(p,q)$-form on $R$. For example by (0.3), 
\begin{equation*}
\phi(\Omega_A)=\phi\left(A^{-1}(z)A_j)^*A^{-1}(z)A_k\right)dz_k\wedge d\overline{z}_j.\tag{1.8}
\end{equation*}
Moreover, $\omega^{p,q}(z)=0$ if and only if $\phi(\omega^{p,q}(z))=0$ for every $\phi\in\bb^*$. We let the set of $\bb$-valued $(p, q)$-forms on $R$ be denoted by $\omega^{p,q}(R,\bb)$.The fact
\[\partial\phi(\omega^{p,q}(z))=\phi(\partial \omega^{p,q}(z))\]
is nicely expressed in the following commuting diagram
\begin{equation*}
\begin{tikzcd}
  \Lambda^{p,q}(R,\bb)\arrow{r}{\partial_p} \arrow{d}{\phi} & \Lambda^{p+1,q}(R,\bb) \arrow{d}{\phi} \\
  \Lambda^{p,q}(R,\C) \arrow{r}{\partial_p} & \Lambda^{p+1,q}(R,\C), 
\end{tikzcd}\tag{1.9}
\end{equation*}
and the parallel diagram regarding $\overline{\partial}$. Here the subscript $p$ in $\partial_p$ is to indicate the space on which $\partial$ acts. Since $d=\partial+\overline{\partial}$, we see that $\omega^{p,q}(z)$ is closed, i.e. $d\omega^{p,q}=0$, if and only if $\phi(\omega^{p,q}(z))$ is closed for every $\phi\in \bb^*$. Since by \cite{HY, ZKKP} the joint resolvent set $P^c(A)$ is a union of domains of holomorphy, we shall only consider the case when $R$ is a domain of holomorphy, in which case its de Rham cohomology $H^p(R,\C)$ (or $H^p(R,\bb)$) can be computed through holomorphic forms. We denote the set of $\bb$-valued holomorphic $p$-forms by $\Lambda^p(R,\bb)$. Note again that in this case $d=\partial$. We refer the readers to \cite{Ra} for details on domains of holomorphy and holomorphic forms. 

Since $\partial_p\partial_{p-1}=0$, one naturally has $\bb$-valued cohomology groups 
\begin{equation*}
H^{p}(R, \bb)= \ker \partial_p/im \partial_{p-1},\ \ p\geq 0,\tag{1.10}
\end{equation*}
where $\partial _{-1}$ is the trivial inclusion map $\{0\}\to \Lambda^0(R,\bb)$. It then follows from the following commuting diagram for holomorphic forms
\begin{equation*}
\begin{tikzcd}
  \Lambda^{p}(R,\bb)\arrow{r}{\partial_p} \arrow{d}{\phi} & \Lambda^{p+1}(R,\bb) \arrow{d}{\phi} \\
  \Lambda^{p}(R,\C) \arrow{r}{\partial_p} & \Lambda^{p+1}(R,\C), 
\end{tikzcd}\tag{1.11}
\end{equation*}
that a form $\omega$ is in $H^{p}(R,\bb)$ if and only if $\phi(\omega)$ is in $\Lambda^{p}(R,\C)$ for every $\phi\in \bb^*$. It was shown in \cite{Ya} that if $\bb$ is a Banach algebra with a trace $\phi$, then for every tuple $A$ such that $P^c(A)$ is nonempty, the $1$-form
$\phi(\omega_A)$ is a nontrivial element in $H^1(P^c(A),\C)$. The above observations thus lead to the following 
\begin{corr}
Let $\bb$ be a unital Banach algebra with a trace and $A$ be a commuting tuple of elements in $\bb$. Then
$\omega_A$ is a nontrivial element in $H^{1}(P^c(A),\bb)$.
\end{corr}
\pf By Proposition 1.1, we see that $\omega_A$ is closed when $A$ is a commuting tuple. If it were exact, then there exists $\bb$-valued smooth function $f(z)$ on $P^c(A)$ such that 
$\omega_A=df$. If $\phi$ is the trace, then it follows that \[\phi(\omega_A)=\phi(df)=d\phi(f).\]
Since $\phi(f)$ is globally defined on $P^c(A)$, the $1$-form $\phi(\omega_A)$ is a trivial element in $H^{1}(P^c(A),\C)$ contradicting with the above mentioned fact in \cite{Ya} that $\phi(\omega_A)$ is a nontrivial element in $H^1(P^c(A),\C)$.
\zb

Proposition 1.1 indicates that the commutativity of tuples $A$ has a natural homological intepretation.

\section{Hermitian metric on $P^c(A)$}

Recall that for a $C^*$-algebra $\bb$, a smooth $\bb$-valued $(1,1)$-form \[\Phi(z)=\Phi_{jk}(z)dz_j\wedge d\overline{z}_k\] on a domain $R$ is said to be Hermitian if for any set of $n$ complex numbers $v_1,\ v_2,\ \dots,\ v_n$ the double sum $v_j\Phi_{jk}(z)\overline{v_k}$ is a positive element in $\bb$ for every $z\in R$. It follows that if $\phi\in \bb^*$ is a positive linear functional, then the matrix $\big(\phi(\Phi_{jk}(z))\big)$ is positive semi-definite for every $z\in R$. By (1.4), we see that the fundamental $(1,1)$-form $\Omega_A=-\omega_A^*\wedge\omega_A$ is Hermitian. This section defines some natural Hermitian metrics on $P^c(A)$ through $\Omega_A$ and its pairing with certain states on $\bb$.

A bounded linear functional $\phi$ on a $C^*$-algebra ${\mathcal B}$ is called a state if it is positive and $\phi(I)=1$. Consider a tuple $A=(A_1, A_2, \cdots, A_n)$ of elements in ${\mathcal B}$. Given a state $\phi$, then on $P^c(A)$ we compute that
\begin{align*}
\phi(\Omega_A(z))&=-\phi (\omega_A^*\wedge \omega_A)\\
&=\phi(A_k^*(A^{-1}(z))^*A^{-1}(z)A_j)dz_j\wedge d\bar{z}_k\\
&:=g_{jk}(z)dz_j\wedge d\bar{z}_k.\tag{2.1}
\end{align*}
As observed in the previous paragraph, the $n\times n$ matrix $g(z)=(g_{jk}(z))$ is positive semi-definite for every $z\in P^c(A)$. So it becomes a natural question for which $\phi$ the matrix $g(z)$ is positive definite on $P^c(A)$. In other words, 
when does $\phi(\Omega_A)$ define a Hermitian metric on $P^c(A)$? To answer this question, we consider the operator space ${\mathcal H}_A=span\{A_1,\ A_2,\ \cdots,\ A_n\}$. A state $\phi$ on $\bb$ is said to be {\em faithful on} ${\mathcal H}_A$ if for every nonzero element $h\in {\mathcal H}_A$ we have $\phi(h^*h)>0$.

\begin{prop}
Let $\phi$ be a state on a $C^*$-algebra $\bb$. Then for any tuple $A$ the $(1, 1)$-form $\phi(\Omega_A)$ induces a Hermitian metric on $P^c(A)$ if and only if $\phi$ is faithful on ${\mathcal H}_A$.
\end{prop}
{\bf Proof}. For every fixed $z\in P^c(A)$, since $A(z)=z_1A_1+z_2A_2+\cdots +z_nA_n$ is invertible, there are numbers $\alpha>\beta>0$ such that
\[\beta I \leq (A^{-1}(z))^*A^{-1}(z) \leq \alpha I.\]
Note that $\alpha$ and $\beta$ depend on $z$, but it is not important here. For every nonzero vector $v=(v_1,v_2,\dots,v_n)\in \C^n$, one checks that
\begin{align*}
v_jg_{jk}(z)\bar{v_k}&=v_j\phi(A_k^*(A^{-1}(z))^*A^{-1}(z)A_j)\bar{v_k}\\
&=\phi(\bar{v_k}A_k^*(A^{-1}(z))^*A^{-1}(z)v_jA_j)\\
&=\phi \left((A^{-1}(z)A(v))^*(A^{-1}(z)A(v)\right)=\phi \left(A^*(v)(A^{-1}(z))^*A^{-1}(z)A(v)\right).
\end{align*}
Since \[\beta A^*(v)A(v)\leq A^*(v)(A^{-1}(z))^*A^{-1}(z)A(v)\leq \alpha A^*(v)A(v)\]
and $\phi$ is positive, it follows that
\[\beta \phi(A^*(v)A(v))\leq v_jg_{jk}(z)\bar{v_k}\leq \alpha \phi(A^*(v)A(v)).\]
Hence $(g_{jk}(z))$ is positive definite for every $z\in P^c(A)$ if and only if $\phi$ is faithful on ${\mathcal H}_A$.
\zb

\begin{exam}
A state $\phi$ on a $C^*$-algebra $\bb$ is called a faithful tracial state if $\phi(ab)=\phi(ba)$ for any $a,\ b\in \bb$ and $\phi(a^*a)>0$ when $a\neq 0$. If $\phi$ is a faithful tracial state, then it is clearly a faithful state on ${\mathcal H}_A$. In this case, we shall write the metric matrix $(g_{ij}(z))$ as $g_A(z)$.
\end{exam}

\begin{exam}
Assume $\bb$ is a $C^*$-algebra of operators acting on a Hilbert space ${\mathcal H}$. Let $\phi_x$ be a vector state defined on $\bb$ by $\phi_x(a)=\la ax,\ x\ra$, where $a\in \bb$ and $x$ is fixed with $\|x\|=1$. If $\{A_1x,\ A_2x,\ \cdots,\ A_nx\}$ is linearly independent, then for any nonzero vector $(c_1,\ c_2,\ \cdots,\ c_n)\in \C^n$, we have 
\[\phi_x((c_iA_i)^*(c_iA_i))=\|c_iA_ix\|^2>0,\]
i.e. $\phi$ is faithful on ${\mathcal H}_A$. 
\end{exam}

It follows from Propositions 1.1, 2.1 and the commuting diagram (1.9) that for a commuting tuple
 $A=(A_1, A_2, \cdots, A_n)$, and every state $\phi$ faithful on ${\mathcal H}_{\hat{A}}$,
\[d\phi(\Omega_{\hat{A}})=\phi(d\Omega_{\hat{A}})=0.\]
Hence $\phi(\Omega_{\hat{A}})$ defines a K\"{a}hler metric on $P^c({\hat{A}})$. Interestingly, the converse is also true if $\phi$ is faithful on the entire $\bb$.

\begin{thm}
Let $\phi$ be a faithful state on $\bb$. Then $\hat{A}=(I, A_1, A_2, \cdots, A_n)$ is commuting if and only if $\phi(\Omega_{\hat{A}})$ is K\"{a}hler.
\end{thm}
{\bf Proof}. It only remains to check sufficiency. If $\phi(\Omega_{\hat{A}})$ is K\"{a}hler, then we have
$0=d\phi(\Omega_{\hat{A}})=\phi(d\Omega_{\hat{A}}).$ Writting $I$ as $A_0$, then by (1.7), we have
\begin{align*}
0&=\phi(\omega^*_{\hat{A}}\wedge \partial \omega_{\hat{A}})\\
&=\sum_{i}\sum_{0\leq j<k}\phi\left(({\hat{A}}^{-1}(z)A_i)^*[{\hat{A}}^{-1}(z)A_j,\ {\hat{A}}^{-1}(z)A_k]\right)d\bar{z}_i\wedge dz_j\wedge dz_k,
\end{align*}
which implies 
\begin{equation*}
\phi\left(({\hat{A}}^{-1}(z)A_i)^*[{\hat{A}}^{-1}(z)A_j,\ {\hat{A}}^{-1}(z)A_k]\right)=0,\ \ \forall 0\leq i, j, k\leq n.\tag{2.2}
\end{equation*}
Applying $\bar{\partial}_m$ to both sides of (2.2) and using the fact $\bar{\partial}_m{\hat{A}}^{-1}(z)=0$ and
\[\bar{\partial}_m ({\hat{A}}^{-1}(z))^*=(\partial_m{\hat{A}}^{-1}(z))^*=-({\hat{A}}^{-1}(z)A_m{\hat{A}}^{-1}(z))^*\]
we have
\begin{equation*}
\phi\left(({\hat{A}}^{-1}(z)A_m{\hat{A}}^{-1}(z)A_i)^*[{\hat{A}}^{-1}(z)A_j,\ {\hat{A}}^{-1}(z)A_k]\right)=0,\ \ \forall 0\leq i, j, k, m\leq n.\tag{2.3}
\end{equation*}
In (2.3), setting $(m, i)$ to $(j, k)$, then to $(k, j)$ and taking the difference, we have
\begin{equation*}
\phi\left([{\hat{A}}^{-1}(z)A_j,\ {\hat{A}}^{-1}(z)A_k]^*[{\hat{A}}^{-1}(z)A_j,\ {\hat{A}}^{-1}(z)A_k]\right)=0,\ \ \forall j, k.
\end{equation*}
Since $\phi$ is faithful, $[{\hat{A}}^{-1}(z)A_j,\ {\hat{A}}^{-1}(z)A_k]=0,  \forall j, k.$
Setting $z$ to $(1, 0, 0, \cdots, 0)$, we have $[A_j, A_k]=0$ for all $j, k$.
\zb

Again, as indicated by Example 1.2 and the remarks before it, the requirement that $I$ be in the tuple $A$ (or in $\mathcal{H}_A$) is indispensible in Theorem 2.4. We end this section by a question motivated by Theorem 2.4.\\

{\bf Problem 1}. Find conditions on commuting tuple $A$ and state $\phi$ such that $\phi(\Omega_A)$ defines a Ricci-flat metric (i.e. Calabi-Yau metric) on $P^c(A)$.

\section{Completeness}

Once a metric is defined on a set, an immediate question is whether the set is complete under the metric. In this section we will study this problem for $P^c(A)$ with respect to the metric defined by $g_A$ as in Examples 2.2. 
Since $\omega_A(z)=A^{-1}(z)dA(z)$ resembles the derivative of logarithmic function, it is natural to expect that $\log$ shall play an important role here.

In \cite{FK}, Fuglede and Kadison defined the following notion of determinant \[\det x=\exp(\phi(\log |x|))\]
for invertible elements $x$ in a finite Von Neumann algebra $ {\mathcal B}$ with a normalized trace $\phi$.
Here $|x|=\sqrt{x^*x} $ and $\log |x|$ is defined by the functional calculus, i.e. 
\[\log |x|=\int_{\sigma(|x|)}\log \lb dE(\lb),\]
where $E(\lb)$ is the associated projection-valued spectral measure. Then
\begin{align}
\phi(\log |x|)=\int_{\sigma(|x|)}\log \lb d\phi(E(\lb)).
\end{align}
Since $x$ is invertible, the spectrum $\sigma(|x|)$ is bounded away from $0$. Hence the integral in (3.1) is greater than $-\infty$, which means $\det x\neq 0$. However, there are non-invertible (singular) elements $x$ for which the improper integral in (3.1) is convergent (hence $\det x\neq 0$). This is essentially due to the absolute convergence of \[\int_0^1\log tdt.\]
An example of such element is given in \cite{FK}. For elements $x$ such that the integral in (3.1) is equal to $-\infty$, $\det x$ is naturally defined to be $0$. This extends $\det$ to all  elements in $\bb$, and by \cite{FK} it is continuous at non-singular elements and upper semi-continuous at singular elements with respect to the norm topology of $\bb$. 

FK-determinant has been well-studied in many papers, and we refer readers to \cite{Ha} for a survey. In particular, it was generalized to $C^*$-algebras in \cite{HS} as follows. 
Assume $GL(\bb)$ is path-connected, $x\in GL(\bb)$ and $x(t)$ is a piece-wise smooth path in $GL(\bb)$ such that $x(0)=I$ and $x(1)=x$, then when the quantity 
\begin{equation}
{\text{det}_{*}} x=\exp \left(\int_0^1 \phi(x^{-1}(t)dx(t))\right) 
\end{equation}
is independent of path $x(t)$, it defines a notion of determinant for invertible elements $x$. One sees that $\text{det}_{*}$ may take on complex values, and it is shown in \cite{Ha} that
\[|{\text{det}_{*}} x|=\exp  \left(Re(\int_0^1 \phi\big(x^{-1}(t)dx(t)\big)\right)=\det x.\]
 Clearly, $\det I=1$, and for a fixed $0\leq s\leq 1$ we have \[ {\text{det}_{*}} x(s)= \exp\big(\int_0^s\phi \left(x^{-1}(t)dx(t)\right)\big),\]
so it follows that
 \begin{equation}
\phi(x^{-1}(t)dx(t))=d \log {\text{det}_{*}} x(t).
\end{equation}
Here $d$ is the differential with respect to $t$. The following definition is needed for our study.\\

{\bf Definition.} Consider a $C^*$-algebra ${\mathcal B}$ with a faithful tracial state $\phi$. Then an element $x$ will be called $\phi$-singular if its FK-determinant $\det x=0$. Furthermore, for a tuple $A$ of elements in ${\mathcal B}$, a point $p\in P(A)$ is said to be $\phi$-singular if $A(p)$ is $\phi$-singular. \\

Now we begin to study the completeness of the metric defined by $g_A$ on $P^c(A)$ in Example 2.2. Let $[P^c(A)]$ denote the completion of $P^c(A)$ with respect to the metric $g_A$. First, consider two points $p$ and $q$ that lie in the same connected component of $P^c(A)$, and let $\gamma=\{z(t):\ 0\leq t\leq 1\}$ be a piecewise smooth path such that $z(0)=p$ and $z(1)=q$. The length of $\gamma$ with respect to the metric $g_A=(g_{ij})$ is
\[L(\gamma)=\int_0^1\sqrt{z'_i(t)g_{ij}(z(t))\overline{z'_j(t)}}dt.\]
The distance from $p$ to $q$ with respect to the metric $g_A$ is the infimum of $L(\gamma)$ over all such paths, i.e.,
\[dist(p, q)=\inf_{\gamma}L(\gamma).\]

\begin{lem}
Let $\bb$ be a $C^*$-algebra with a faithful tracial state $\phi$. If $p$ and $q$ are in the same connected component of $P^c(A)$ then \[dist(p,\ q)\geq\mid \phi( \log |A(p)|)-\phi(\log |A(q)|)\mid.\]
\end{lem}

{\bf Proof}. First of all, by (2.1) we have $g_{ij}(z)=\phi\big((A^{-1}(z)A_j)^*A^{-1}(z)A_i\big)$. Let $p$ and $q$ be two points in the connected compnent $U$ of $P^c(A)$,
 and let $\gamma=\{z(t)|\ 0\leq t\leq 1\}$ be a piecewise smooth path such that $z(0)=p$ and $z(1)=q$. Then on $\gamma$ we have $A'(z(t))=z_j'(t)A_j$, and the length of $\gamma$ can be computed as
\begin{align*}
L(\gamma)&=\int_0^1\sqrt{z'_i(t)g_{ij}(z(t))\overline{z'_j(t)}}dt\\
&=\int_0^1\sqrt{z'_i(t) \phi\big((A^{-1}(z)A_j)^*A^{-1}(z)A_i\big)  \overline{z'_j(t)}}dt\\
&=\int_0^1\sqrt{\phi\big((A^{-1}(z(t))z'_j(t)A_j)^*A^{-1}(z(t))z'_i(t)A_i\big)}dt\\
&=\int_0^1\sqrt{\phi\big(\left(A^{-1}(z(t))A'(z(t)\right)^*A^{-1}(z(t))A'(z(t))\big)}dt. \tag{3.4}
\end{align*}
Since $\phi$ is a state, for every $a,\ b\in \bb$ one has $|\phi(ab)|^2\leq \phi(a^*a)\phi(b^*b)$. In particular, $|\phi(a)|^2\leq \phi(a^*a)$.
Hence 
\begin{align*}
L(\gamma)&\geq \int_0^1|\phi (A^{-1}(z(t))A'(z(t))|dt=\int_0^1|\phi [A^{-1}(z(t))dA(z(t))]|\\
&=\int_0^1|\phi (\omega_A(z(t))|\geq |\int_0^1\phi (\omega_A(z(t))|.
\end{align*}
By (3.3) we have
\begin{align*}
L(\gamma)&\geq |\int_0^1d\log {\text{det}_{*}} A(z(t))|\\
&=\mid \log \text{det}_{*} A(p)-\log \text{det}_{*} A(q)\mid\\
&=\mid \big(\log |\text{det}_{*} A(p)|-\log |\text{det}_{*} A(q)|\big)+i\big(Arg \text{det}_{*} A(p)-Arg \text{det}_{*} A(q)\big)\mid\\
&\geq \mid \log |\text{det}_{*} A(p)|-\log |\text{det}_{*} A(q)|\mid=\mid \log \det A(p)-\log{\det} A(q)\mid.\tag{3.5}
\end{align*}
Since $\gamma$ is arbitary, we have 
\begin{align*}
dist(p,\ q)=\inf_{\gamma}L(\gamma)& \geq \mid \log {\det} A(p)-\log {\det} A(q)\mid\\
&=\mid \phi( \log |A(p)|)-\phi(\log |A(q)|)\mid. \tag{3.6}
\end{align*}
\zb

The proof of Lemma 3.1 can be modified a little to accommodate the case in which $q$ is a boundary point of a path-connected component of $P^c(A)$. First, we observe that in this case $q\in P(A)$. Suppose there exists a piece-wise smooth path $z(t)$ such that $z(t)\in P^c(A),\ 0\leq t<1$ but $q=z(1)\in P(A)$. Then by the proof above we have
\begin{align*}
L(\gamma)&\geq \lim_{s\to 1^{-}}\int_0^s|\phi(\omega_A(z(t))|\\
&\geq \limsup_{s\to 1^{-}}\mid \log {\det} A(p)-\log{\det} A(z(s))\mid.\tag{3.7}
\end{align*}
If $q$ is a $\phi$-singular point, we have $\det A(z(1))=\det A(q)=0$. Note that $0$ is the minimum possible value for $\det$. Hence the upper semi-continuity of $\det A(z(s))$ implies \[\lim_{s\to 1}\det A(z(s))=0.\] Therefore by (3.7) we have $L(\gamma)=\infty$ for every such path $\gamma$ in $P^c(A)$ connecting $p\in P^c(A)$ to $q\in P(A)$, and consequently  
 \[dist(p,\ q)=\inf_{\gamma}L(\gamma)=\infty.\]
Since $q$ has infinite distance to every $p\in P^c(A)$, we have the following
\begin{thm}
For a tuple $A$ in a $C^*$-algebra $\bb$ with a faithful tracial state $\phi$, if $q\in \partial P^c(A)$ is $\phi$-singular, then $q\notin [P^c(A)]$.
\end{thm}

If $\bb$ is a matrix subalgebra in $M_k(\C)$, we let $Tr$ and $det$ stand for the ordinary trace and respectively determinant of $k\times k$ matrices. Then for any $n$-tuple $A$, the joint spectrum $P(A)$ is the hypersurface $\{det A(z)=0\}$ in $\C^n$. Clearly, $\partial P^c(A)=P(A)$. Let $\phi=\frac{1}{k}Tr$. Then $\phi$ is a tracial state on $M_k(\C)$, and in this case the FK-determinant and the usual determinant satisfies the relation ${\det}x=|det x|^{1/k}$ (cf. \cite{FK}). So in this case every point in $P(A)$ is $\phi$-singular. Since the factor $1/k$ is not important, we have the following

\begin{corr}
For every tuple $A$ of $k\times k$ matrices, the metric on $P^c(A)$ defined by $Tr(\Omega_A)$ is complete.
\end{corr}

The general linear group $GL(k, \C)$ is an open subset in $M_k(\C)$, hence it is clearly not complete with respect to the usual Euclidean metric. Corollary 3.3 implies that we can endow $GL(k)$ with a complete metric.

\begin{corr}
Let $A_1,\ A_2,\ \cdots,\ A_n$ be a basis for the vector space $M_k(\C)$, where $n=k^2$. Then $Tr(\Omega_A)$ defines a complete, left invariant, Ricci flat, but non-K\"{a}hler metric on the complex general linear group $GL(k)$.
\end{corr}
{\bf Proof}. First, since $A_1,\ A_2,\ \cdots,\ A_n$ is a basis for the matrix algebra $M_k(\C)$, every matrix is of the form $A(z)$ for some $z$. Hence $z\in P^c(A)$ if and only if $A(z)\in GL(k)$. And this provides an identification of $P^c(A)$ with $GL(k)$. So in particular $Tr(\Omega_A)$ defines a metric on $GL(k)$. Recall that 
\[Tr(\Omega_A)=-Tr(\omega_A^*\wedge \omega_A)=Tr\big((A^{-1}(z)A_j)^*A^{-1}(z)A_i\big)dz_i\wedge d\bar{z}_j.\]
For $1\leq i\leq n$, we write $A_i=(a^i_1, a^i_2, \cdots, a^i_k)$, where $a^i_m,\ 1\leq m\leq k,$ is the $m$th column of $A_i$. Then
\begin{align*}
&(A^{-1}(z)A_j)^*A^{-1}(z)A_i\\
&=\big(A^{-1}(z)a_1^j, A^{-1}(z)a_2^j, \dots, A^{-1}(z)a_k^j\big)^* \big(A^{-1}(z)a_1^i, A^{-1}(z)a_2^i, \dots, A^{-1}(z)a_k^i\big),
\end{align*}
and hence
\begin{equation}
Tr\big((A^{-1}(z)A_j)^*A^{-1}(z)A_i\big)=\sum_m\la A^{-1}(z)a_m^i,\ A^{-1}(z)a_m^j\ra_{\C^k}.
\end{equation}

Let $\alpha$ be the $n\times n$-matrix $\left(a_m^i\right)_{1\leq m \leq k,1\leq i \leq n}$, whose $i$-th column (in $\C^n$) is denoted by $\alpha_i$, and let $I_k$ be the $k\times k$ identity matrix. Then $A^{-1}(z)\otimes I_k$ is the $n\times n$ block matrix with $A^{-1}(z)$ on the diagonal. So we can write the sum in (3.4) as 
\begin{equation}
\la(A^{-1}(z)\otimes I_k)\alpha_i,\ (A^{-1}(z)\otimes I_k)\alpha_j\ra_{\C^n}.
\end{equation}

Since $A_1,\ A_2,\ \cdots,\ A_n$ are linearly independent, $\alpha$ is invertible, and one checks by direct computation that the metric matrix
\begin{align*}
g(z)&=\left( Tr((A^{-1}(z)A_j)^*A^{-1}(z)A_i\right)_{n\times n}\\
&=\left(\la(A^{-1}(z)\otimes I_k)\alpha_i,\ (A^{-1}(z)\otimes I_k)\alpha_j\ra\right)_{n\times n}\\
&=\alpha^* (A^{-1}(z)\otimes I_k)^* (A^{-1}(z)\otimes I_k)\alpha.
\end{align*}
It follows that $det g(z)=|det \alpha|^2|det A^{-k}(z)|^{2}$. Since $det A^{-k}(z)$ is holomorphic on $P^c(A)$, the metric defined by $Tr\Omega_A$ is Ricci flat by the remark after (1.1).

Moreover, because $A_1,\ A_2,\ \cdots,\ A_n$ is a basis for $M_k(\C)$, for every fixed $x\in GL(k)$, $xA(z)=A(w)$, where $w\in P^c(A)$ and is uniquely determined by $z$.  This gives rise to a natural action $L$ of $GL(k)$ on $P^c(A)$ given by $L_x(z)=w$. Since by (0.2), $\omega_A$ is invariant under left action of $GL(k)$, the metric defined by $Tr\Omega_A$ is invariant under the action of $L$. The completeness of the metric follows from Corollary 3.3. Since $I$ is of the form $A(z)$ for some $z\in \C^n$ in this case, the non-K\"{a}hlerness follows from Theorem 2.4.
\zb \\

We end the paper by posting\\

{\bf Problem 2.} Is the converse of Theorem 3.2 true?

\subsection*{Acknowledgment}
The second-named author would like to thank Guoliang Yu and the Department of Mathematics at Texas A\&M University for their support and hospitality during his visit.

\vspace{5mm}

\end{document}